\documentclass[12pt]{amsart}
\usepackage{amsthm}
\usepackage{amsmath}
\usepackage{amsfonts}
\usepackage{amssymb}
\usepackage{amscd}

\theoremstyle{definition}
\newtheorem{Thm}{Theorem}

\newtheorem{Prop}[Thm]{Proposition}
\newtheorem*{Rem}{Remark}

\newtheorem*{T7}{Theorem 1}

\newtheorem{Def}[Thm]{Definition}

\newcommand{\delbar}{\bar\partial}
\newcommand{\Sum}{\sum}

\newcommand{\Dim}{\operatorname{dim}}

\newcommand{\tensor}{\otimes}
\newcommand{\st}{\mid}

\newcommand{\Laplace}{\Delta}
\newcommand{\by}{\times}

\newcommand{\vardelta}{\delta}
\newcommand{\cG}{\mathcal G}

\newcommand{\End}{\operatorname{End}}

\newcommand{\tr}{\operatorname{tr}}

\newcommand{\R}{{\mathbb R}}
\newcommand{\bC}{{\mathbb C}}
\newcommand{\bZ}{{\mathbb Z}}
\newcommand{\bR}{{\mathbb R}}
\newcommand{\bH}{{\mathbb H}}

\newcommand{\too}{\longrightarrow}

\newcommand{\func}{\operatorname}

\begin{document}

\title[Analytic Torsion]{Analytic Torsion for Quaternionic Manifolds and Related Topics}
\date{September 1997}
\author[N.C.Leung]{Naichung Conan Leung}
\address{Department of Mathematics\\
University of Minnesota\\
Minneapolis, Minnesota 55455}
\email{leung@math.umn.edu}
\author[S. Yi]{Sangkug Yi}
\email{yi@math.umn.edu}
\begin{abstract}
In this paper we show that the Ray-Singer complex analytic torsion is
trivial for even dimensional Calabi-Yau manifolds. 

Then we define
the quaternionic analytic torsion for quaternionic manifolds and prove that they
are metric independent. 

In dimension four, the quaternionic analytic torsion equals
to the self-dual analytic torsion. For higher dimensional manifolds, the 
self-dual analytic torsion is a conformal invariant. 
\end{abstract}
\maketitle

\allowdisplaybreaks
\section{Introduction}

In the seminal papers $\left( \left[ RS1\right] ,\left[ RS2\right] \right) ,$
Ray and Singer defined analytic torsions for real and complex manifolds as
determinants of the deRham complex and Dolbeault complex respectively.
As an analog to the vanishing of real analytic torsions for even dimensional
manifolds, we prove the following vanishing result for complex analytic
torsions.

\begin{Thm}
If $M$ is an even dimensional Calabi-Yau manifold, then 
\[
\log \tau _{\bC}\left( M,V\right) =0,
\]
where $\tau _{\bC}\left( M,V\right) $ is the complex analytic torsion
for $M$ with coefficient in the unitary flat bundle $V.$
\end{Thm}

In particular the complex analytic torsion is always trivial for HyperK\"ahler 
manifolds. Naturally we look for refinement of the complex analytic 
torsion. In section \ref{qtorsion},  we define a
quaternionic analytic torsion $\tau _{\bH}\left( M,V\right) $ for
quaternionic manifolds. 

Recall that a $4n$-dimensional manifold $M$ is
called a quaternionic manifold if its tangent bundle admits a torsion-free $%
GL\left( n,\bH\right) Sp(1)$ connection. Let $E$ and $H$ denote the 
associated bundles to the frame bundle of $M$ with respect to standard
representations of $GL(n,\bH)$  and $Sp(1)$.

\begin{Thm}
Let $V$ be a unitary flat bundle over a quaternionic
manifold $M$.  For any compatible metric on $M,$ we
define $\tau _{\bH}(M,V) $ to be the regularized
determinant of the following elliptic complex: 
\[
0\rightarrow A^{0}\stackrel{D}{\longrightarrow }A^{1}\stackrel{D}{\longrightarrow }%
\dots\stackrel{D}{\longrightarrow }A^{2n}\rightarrow 0,
\]
where $A^{k}=\Gamma(M,\Lambda ^{k}E\otimes S^{k}H\otimes V)$.

Then the ratio $\tau _{\bH}\left( M,V_{1}\right) /\tau _{\bH}\left(
M,V_{2}\right) $ is independent of the choice of the compatible metric on $M$
provided that their corresponding complexes have trivial cohomology groups.
\end{Thm}

The above complex was introduced by S. M. Salamon (\cite{S1}).

\begin{Rem}
The assumption about trivial cohomology groups is not essential. Otherwise
we can write down the variation of $\tau_\bH$ 
in terms of the change of volume forms on the determinant of cohomology groups.
\end{Rem}

Notice that every Riemannian metric on a four manifold $M$ determines a
quaternionic structure on $M$. 
The above elliptic complex is the same as the self-dual complex: 
\[
0\rightarrow \Omega ^{0}\left( M,V\right) \stackrel{d}{\longrightarrow }\Omega
^{1}\left( M,V\right) \stackrel{\sqrt{2}P_{+}d}{\longrightarrow }\Omega
_{+}^{2}\left( M,V\right) \rightarrow 0.
\]

Because conformally equivalent metrics determine the same quaternionic
structure on $M,$ above theorem shows that the regularized determinant of the
self-dual complex of a four manifold depends only on the conformal class of
the metric. More generally, the following is true:

\begin{Thm}
We assume that $M$ is $4n$-dimensional closed Riemannian manifold and $V$ 
is an orthogonal flat vector bundle over $M$. Let 
$$\tau_{\text{SD}}( M,V) = \sum\limits_{q=0}^{2n-1}\left(
-1\right) ^{q+1}q\log \det \left(\bigtriangleup _q\right) -n\log
\det \left(\bigtriangleup_{2n}\right). $$ 
Then 

(i) $\tau_{\text{SD}}(M,V)$ is the regularized determinant of the
following self-dual complex: 
\begin{align*}
0\rightarrow \Omega ^{0}(M,V) \stackrel{d}{\rightarrow }
\dots \stackrel{d}\rightarrow\Omega^{2n-1}(M,V)\stackrel{\sqrt{2}P_{+}d}{\rightarrow}\Omega_{+}^{2n}(M,V) \stackrel{}{\rightarrow }
0.
\end{align*}

(ii) $\tau _{\text{SD}}\left( M,V_{1}\right) /\tau _{\text{SD}}\left( M,V_{2}\right) $
depends only on the conformal class of the Riemannian metric provided that
the above self-dual complex has trivial cohomology groups.
\end{Thm}

A similar result on conformal invariance of self-dual analytic torsion has been known to J. Cheeger for many years. After this paper is finished, the author is informed by Bismut that he also noticed the result of theorem one before by  using Serre dualit
y for Quillen metrics.

%
%
\section{Vanishing of  Complex Analytic Torsions}

Ray and Singer \cite{RS1} defined the real analytic torsion $\tau_\bR(M)$ for a compact oriented
Riemannian manifold $M$ of dimension $n$, using zeta functions of Laplacian
and showed that they are metric independent.
Moreover, $\tau_\bR$ is trivial when $M$ is an even dimensional manifold. In
this section, we shall prove a complex analog of their theorem.
To begin, we first recall their construction of the real analytic torsion,
$$ \tau_\bR(M,V) = \exp\left(\frac12\sum^n_{q=0}(-1)^qq\zeta'_{q,V}(0)\right), $$
where $V$ denotes both an orthogonal representation of $\pi_1(M)$ and
the flat bundle associated with that representation (this convention will be 
used for the whole paper). Here, $\zeta$ is the zeta
function for the Laplacian $\Delta_q=(d\delta+\delta d)$ defined on the space of $q$-forms
with values in $V$, where $d$ is the
exterior differentiation and $\delta$ is the formal adjoint of $d$.
The zeta function of a Laplacian is defined as follows.
If $\lambda_1\le \lambda_2 \le \dots$ are positive eigenvalues 
of $\Delta_q$ listed in the increasing order counting multiplicities.
Then the series
$$\zeta_q(s)=\sum_{i> 0} \frac1{\lambda_i}=\frac1{\Gamma(s)}
\int_0^\infty t^{s-1}e^{-t\Delta_q}dt$$
is defined for $s$ with ${\func{Re}} s > 1$ provided that $H^q(M,V)=0$. This function of $s$ has the analytic
continuation to the entire complex plane with a continuous derivative
at $s=0$.
Equivalently,
$$\tau_\bR(M,V) = \prod_{0\le q \le n}(\det\thinspace\negthinspace'(\Delta_q))^{(-1)^{q+1}q},$$
where $\det'(\Delta_q)$ is the regularized determinant of $\Delta_q$ defined
as
$\det'(-\Delta_q)=\exp(\zeta'_q(0))$.

We have the theorem of Ray and Singer :

\begin{Thm}\cite{RS1}
Let $M$ be a closed manifold and $V$ be an
orthogonal representation of $\pi_1(M)$ with the 
property that the cohomology of $M$ with coefficients in 
$V$ is trivial. 

Then $\tau_\bR(M,V)$ is metric independent.
\end{Thm}
Let $V^*$ be the dual of $V$ then if we consider the dual complex defined by
the formal adjoint of the exterior differentiation, we have
$$\tau_\bR(M,V)=\tau_\bR(M,V^*)^{(-1)^{(n+1)}}$$
which implies the vanishing result of Ray and Singer when the dimension of $M$
is even. 
\begin{Thm}
When the dimension of $M$ is even, $\tau_\bR(M,V)= 1$, or equivalently,
$\log\tau_\bR(M,V)= 0$.
\end{Thm}

In particular, the real analytic torsion $\tau_\R(M,V)$ vanishes for complex
manifolds. But for a complex manifold the deRham complex has a canonical 
subcomplex, namely, the Dolbeault complex. In the sequel \cite{RS2} to their
paper about the real analytic torsion, Ray and Singer defined complex analytic torsion
$\tau_\bC(M,V)$ of a $n$-dimensional complex manifold $M$ as the regularized determinant of the Dolbeault complex
$$0\too \Omega^{0,0}(M,V) \stackrel{\delbar}\too\Omega^{0,1}(M,V)\stackrel{\delbar}
\too\dotsi\stackrel{\delbar}\too\Omega^{0,n}(M,V)\too 0,$$
using the $\delbar$-Laplacian $\Delta_{\delbar}=\delbar\delbar^* +\delbar^*\delbar$ 
and its zeta function. Here $V$ is the vector bundle associated with a unitary 
representation of $\pi_1(M)$.
The complex analytic torsion is also known as the holomorphic
torsion. Ray and Singer showed that the ratio between two complex analytic torsions
corresponding to different unitary representations of the fundamental group is 
independent of the Hermitian metric chosen, if all the cohomology of the Dolbeault complexes vanish.

As we have seen in the real case, a certain duality property for the de Rham
complex gives the vanishing of the real analytic torsion when the dimension of the
manifold is even. The dual map in the de Rham complex is given by the Hodge-$*$
 operator which is defined using the volume form. In the complex manifold
case, we can define a similar duality map if there exists a parallel holomorphic
volume form. Hence this idea leads us to the condition of Calabi-Yau manifolds.

Recall that a compact K\"ahler manifold $M$ is called a Calabi-Yau manifold
if it has trivial canonical line bundle, $K_M={\mathcal{O}}_M$. By the theorem of 
Yau \cite{Y}, such a manifold admits a unique Ricci flat K\"ahler metric in
each K\"ahler class. Equivalently, such a manifold has $SU(n)$ holonomy
where $n=\dim_\bC M$. Using Bochner arguments, we can show that the unique(up to scalar multiplications) holomorphic $n$-form $\Omega$ on $M$ is parallel.
Therefore, $\Omega\wedge\bar\Omega$ is a constant multiple of the volume
form on $M$ and we call $\Omega$ a holomorphic volume form on $M$.

Now we prove a complex analog of the Ray-Singer vanishing theorem for 
real analytic torsions:

\begin{T7}
If $M$ is an even dimensional Calabi-Yau manifold, then
\[
\log \tau _{\bC}\left( M,V\right) =0,
\]
where $\tau _{\bC}\left( M,V\right) $ is the complex analytic torsion
for $M$ with coefficient in the unitary flat bundle $V.$
\end{T7}

\begin{proof}
Let $M$ be a Calabi-Yau manifold with complex dimension $2n$.
Let $H^\lambda_q=\{\phi \in \Omega^{0,q}(M,V) | \Delta_{\bar\partial}\phi = \lambda \phi \}$,
$E^\lambda_q=\{\phi \in H^\lambda_q | \delbar\phi = 0 \}$,
${E'}^\lambda_q=\{\phi \in H^\lambda_q | \delbar^*\phi = 0 \}$.

Then $H^\lambda_q$ is an orthogonal direct sum of $E^\lambda_q$ and ${E'}^\lambda_q$, $H^\lambda_q=E^\lambda_q\oplus {E'}^\lambda_q$.

Let $N^\lambda_q=\Dim E^\lambda_q$ and ${N'}^\lambda_q=\Dim {E'}^\lambda_q$.  
By the isomorphism $\frac1{\sqrt\lambda}\delbar : {E'}^\lambda_q \too E^\lambda_
{q+1}$(with the inverse $\frac1{\sqrt\lambda}\delbar^*)$,

$$\Dim H^\lambda_q=N^\lambda_q+{N'}^\lambda_q=N^\lambda_q+N^\lambda_{q+1}={N'}^\lambda_q+{N'}^\lambda_{q-1}.$$

Hence,
\begin{align*}
\zeta_q(s)=&\Sum_{\lambda\neq0}\frac1{\lambda^s}(N^\lambda_q+{N'}^\lambda_{q+1})
=\Sum_{\lambda\neq0}\frac1{\lambda^s}(N^\lambda_q+N^\lambda_{q+1})\\
=&\Sum_{\lambda\neq0}\frac1{\lambda^s}({N'}^\lambda_q+{N'}^\lambda_{q-1}).
\end{align*}
\negthinspace and
\begin{align*}
\Sum_{q=0}^{2n}(-1)^qq\zeta_q(s)=\Sum_{q=1}^{2n}(-1)^q\Sum_{\lambda>0}\frac1{\lambda^s}N^\lambda_q=
\Sum_{q=0}^{2n-1}(-1)^{q+1}\Sum_{\lambda>0}\frac1{\lambda^s}{N'}^\lambda_q.
\end{align*}

Let $\Omega$ be the holomorphic volume form on $M$. Then we have an isomorphism
$A:\Omega^{0,q}\too\Omega^{2n,q}$ given by wedging with $\Omega$.
Since $V$ is hermitian, the dual Dolbeault complex is identified with the
original complex. Hence we get the isomorphism
$A^{-1}\circ\bar*:E^\lambda_q\too {E'}^\lambda_{2n-q}$.
Hence we have $N^\lambda_q=N^\lambda_{2n-q}$.
Therefore,
\begin{alignat*}{3}
\Sum_{q=1}^{2n}(-1)^q\Sum_{\lambda>0}\frac1{\lambda^s}N^\lambda_q
&=\Sum_{q=0}^{2n-1}(-1)^{2n-q}\Sum_{\lambda>0}\frac1{\lambda^s}N^\lambda_{2n-q}
&=\Sum_{q=0}^{2n-1}(-1)^{q}\Sum_{\lambda>0}\frac1{\lambda^s}{N'}^\lambda_{q}\\
&=-\Sum_{q=0}^{2n-1}(-1)^{q+1}\Sum_{\lambda>0}\frac1{\lambda^s}{N'}^\lambda_q 
&=-\Sum_{q=1}^{2n}(-1)^q\Sum_{\lambda>0}\frac1{\lambda^s}N^\lambda_q.
\end{alignat*}
Hence the summation is zero and the log of the complex analytic
torsion given by 
$$\log\tau_\bC(M,V)=\Sum_{q=0}^{2n}(-1)^qq\zeta_q'(s)\st_{s=0}$$ vanishes.

\end{proof}

In particular, if the manifold has a holonomy contained in $Sp(n)$, (namely
a HyperK\"ahler manifold) then the manifold has a trivial complex analytic
torsion. In the next section, we will define the quaternionic analytic
torsion for such manifolds. More generally, they are defined for manifolds
with its holonomy contained in $GL(n,\bH)Sp(1)$, namely, quaternionic
manifolds.

\section{Quaternionic Analytic Torsions}\label{qtorsion}
\subsection{Basic Facts about Quaternionic Manifolds}

A 4n-dimensional manifold $M$ is called an almost quaternionic manifold if there is a 
rank 3 subbundle ${\mathcal G}$ of $\End(TM)$ such that for each $x\in M$ there is 
a neighbourhood $U$ of $x$ over which ${\mathcal G}\arrowvert_U$ has a basis $\{I,J,K\}$ of 
almost complex 
structures with $K=IJ=-JI$. Note that this is only a local basis. 
An example of this basis over an open subset of ${\mathbb R}^{4n}$ is, with
the usual identification of tangent space of $\bR^{4n}$ with the ${\mathbb R}^{4n}$ 
itself, 
{\tiny
$$I=\left(\begin{array}{cccc} 0 & -1_n & 0 & 0 \\
		1_n & 0 & 0 & 0 \\
		0 & 0 & 0 & -1_n \\
		0 & 0 & 1_n & 0 \end{array}\right),
J=\left(\begin{array}{cccc} 0 & 0 &-1_n & 0 \\
                0 & 0 & 0 & 1_n  \\
                1_n & 0 & 0 & 0  \\
                0 & -1_n & 0 & 0\end{array}\right), 
K=\left(\begin{array}{cccc} 0 & 0 & 0 & -1_n \\
                0 & 0 & - 1_n & 0 \\
                0 & 1_n & 0 & 0 \\
                1_n & 0 & 0 & 0\end{array}\right)$$ }
\negthinspace\negthinspace where $1_n$ represents the $n\by n$ identity matrix.
In other words,
\begin{Def}
A real 4n-dimensional manifold is almost quaternionic if the structure
group of its principal frame bundle can be reduced from $GL(4n,\bR)$ to $G=GL(n,\bH\,)Sp(1)$.
\end{Def}

 If $M$ is 
a Riemannian manifold, the metric $g$ is said to be compatible with the 
almost quaternionic structure ${\mathcal G}$ if 
$$g(AX,AY)=g(X,Y)$$ for $X,Y\in T_xM_\bC$ and $A \in {\cG}_{x}$
such that $A^2=-1$. Locally, we can  construct a compatible metric $g$ from any 
Riemannian metric $g'$ by defining
$$g(X,Y)=\frac14(g'(X,Y)+g'(IX,IY)+g'(JX,JY)+g'(KX,KY)).$$
for $X,Y \in TM_{\mathbb C}$. Actually, this definition doesn't depend on the
choice of  the basis $\{I,J,K\}$. hence $g$ is defined globally.

Note that the space of all the metric compatible with the given quaternionic
structure is 
contractible. From this fact, when we consider a variation of compatible metric, it is 
enough to look at only local variations. 

Next, we discuss quaternionic manifolds.
\begin{Def}
A real $4n$-dimensional almost quaterionic manifold $M$ is quaternionic
if the principal frame bundle of $M$ admits a torsion-free $GL(n,\bH)Sp(1)$-connection.
\end{Def}

Note that the existence of a torsion-free $G$-connection is only a partial
obstruction to the existence of an integrable $G$-structure.

If $M$ is a quaternionic manifold of dimension $4n$, then its holonomy group
is a subgroup of $GL(n,\bH\,)Sp(1)$ and the reduced frame bundle 
$F$ consists of 
frames $u:{\bH}^n\too T_xM$ which are compatible with the quaternionic 
structure. Locally, $F$ can be lifted to a principal $GL(n,\bH\,)\times 
Sp(1)$-bundle $\tilde F$ which double-covers $F$. This bundle exists globally if 
$(-1,-1)\in GL(n, \bH\,)\times Sp(1)$ acts as the identity.
The obstruction to the global existence of the double cover $\tilde F$ is the
vanishing of the cohomology class $w(M) \in H^2(M,{\mathbb Z}_2)$ 
defined as follows.

Let $[F]$ denote the element of the cohomology group 
$H^1(M;GL(n,{\bH}\,)Sp(1))$ corresponding to the principal 
$GL(n,{\bH}\,)Sp(1)$-bundle $F$ over $M$. 

>From the short exact sequence
$$0\too {\mathbb Z}_2\too GL(n,{\bH}\,)\times 
Sp(1)\too GL(n,{\bH}\,)Sp(1)\too 0$$
we get a coboundary homomorphism
$$\delta:H^1(M;GL(n,{\bH}\,)Sp(1))\too H^2(M;{\mathbb Z}_2).$$
If we define $w(M)=\delta[F] \in H^2(M;{\mathbb Z}_2)$, then the class $w(M)$ is the obstruction to lifting $F$ to 
a principal $(GL(n,{\bH}\,)\times Sp(1))$-bundle $\tilde F$.  

>From the lifting to a principal 
$GL(n,{\bH}\,)\times Sp(1)$ bundle, we get the decomposition
 $TM_{\mathbb C}=TM\tensor_{\mathbb R}{\mathbb C}
=E\tensor H$, where $E$ is the vector bundle associated to the standard
representation of $GL(n,{\bH}\,)$ on $\bH^n=\bC^{2n}$, and $H$ is the vector bundle
associated to the standard representation of $Sp(1)$ on $\bH^n=\bC^{2n}$. 

For $q\le 2n$, we denote 
$$A^q=\bigwedge^qE\tensor S^qH$$
, where $S^qH$ denotes $q$-th
symmetric power of $H$ and define $D$ to be the differential operator defined
as $d$ followed by the projection onto $A^q$. Then by the work of S. M. 
Salamon \cite{S1}, we get an elliptic complex
$$0\too A^0(M)\stackrel{D}\too A^1(M) \stackrel{D}\too \dotsi\stackrel{D}
\too A^{2n}(M) \too 0.$$
\begin{Rem}
 If the class $w(M)=0$, then we have two vector bundles $E$ and $H$. However, 
the existence of $A^q(M)=\bigwedge^qE\tensor S^qH$ and the elliptic complex
$$0\too A^0(M)\stackrel{D}\too A^1(M) \stackrel{D}\too \dotsi\stackrel{D}
\too A^{2n}(M) \too 0,$$
do not depend on $w(M)$ being zero or not. This complex will be called the
quaternionic complex. We define the cohomology group
for this complex $H^q_A(M)=H^q(A^*(M),D)$ for any quaternionic manifold.
\end{Rem}

\begin{Rem}(The four dimension case)
Notice that every Riemannian four manifold $M$ determines a quaternionic
manifold because of $SO\left( 4\right) =Sp\left( 1\right) Sp\left(
1\right) \subset GL\left( 1,\bH\right)Sp(1) $.
Moreover, $w(M)$ equals to the second Stiefel-Whitney class of $M$, $w_2(M)$. Therefore
$w(M)=0$ is equivalent to $w_{2}\left( M\right) =0,$ namely $M$
is a Spin manifold. This identification follows from the following commutative diagram:

\begin{tabular}{ccccc}
$\bZ_2$ & = &$\bZ_2$ & = & $\bZ_2$ \\
$\downarrow$ & & $\downarrow$ & & $\downarrow$ \\
$Spin(4)$ & = & $Sp(1)\times Sp(1)$ & $\subset$ & $GL(1,\bH) \times Sp(1)$ \\
$\downarrow$ & & $\downarrow$ & & $\downarrow$ \\
$SO(4)$ & = & $Sp(1)Sp(1)$ & $\subset$ & $GL(1,\bH)Sp(1)$
\end{tabular}

Furthermore, we have $E=S^{-}$ and $H=S^{+},$ where $S^{\pm }
$ are the positive/negative spinor bundles over $M$ and the
quaternionic complex is the same as the self-dual complex: 
\[
0\rightarrow \Omega ^{0}(M,V) \stackrel{d}{\rightarrow }\Omega
^{1}(M) \stackrel{P_{+}d}{\rightarrow }\Omega _{+}^{2}(
M) \rightarrow 0.
\]
This is because of the canonical identifications $T_{M}^{*}\otimes_{\bR
}\bC=S^{-}\otimes S^{+}$, $\Lambda _{+}^{2}T_{M}^{*}={\mathfrak{su}}
(S^{+}) $ and the fact that $S^{\pm }$ are $SU(2)$
bundles.

In dimension four, people also call an Einstein manifold with quaternionic 
structure a quaternionic manifold.
\end{Rem}

\begin{Rem}
There are two important classes of quaternionic manifolds which had been
extensively studied in the literatures. Namely, they are the quaternionic
K\"ahler manifolds (holonomy inside $Sp(n)Sp(1)$) and the HyperK\"ahler 
manifolds (holonomy inside $Sp(n)$).
\end{Rem}

\subsection{The Definition of Quaternionic Analytic Torsion}

In this section, let the manifold $M$ be a compact quaternionic manifold. 
In other words, we have the tensor product decomposition 
$TM_{\mathbb C}=E\tensor H$ of the complexified tangent bundle of $M$.

If we define the Young diagram $\lambda_q$ for $q>2n$ to be the diagram 
obtained by joining $2n\times 1$ diagram with $\lambda_{q-2n}$,
we can also define $A^q$ for $q$ bigger than $2n$. Using these, we can define
another complex
$$0\longrightarrow A^{2n} \stackrel{d}\longrightarrow \dots \stackrel{d}
\longrightarrow A^{4n} \longrightarrow 0.$$
Consider the pairing,
$$A^q\tensor A^{4n-q}\stackrel{\alpha}\too \bigwedge^{4n}TM_{\mathbb C}=S_{\lambda_{4n}}E
\tensor S_{\lambda_{4n}'}H.$$
Since we have the volume form $vol$ on $\bigwedge^{4n}TM_{\mathbb C}$ and 
the inner product on $A^q$ given by g, we can define 
$\gamma:S_{\lambda_q}\too S_{\lambda_{4n-q}}$ by
$$\alpha(a\tensor \gamma b) = <a,b>vol$$
for $a,b \in A^q$.
We can represent the pairing using Young diagrams.
\vskip.2in
\begin{tabular}{|c|} \hline
	{\tiny 1} \\ \hline
	{\tiny 2}\\ \hline
	:\\ 
	:\\ \hline
        {\tiny q}\\ \hline \end{tabular} \quad$\tensor$
\quad \begin{tabular}{|c|c|c|c|} \hline
	{\tiny 1} & {\tiny 2} & $\dotsi$ & {\tiny q} \\ \hline
	\end{tabular} \quad 
$\tensor$
\quad \begin{tabular}{|c|c|} \hline
	{\tiny 1} & {\tiny 2n+1} \\ \hline
	{\tiny 2} & {\tiny 2n+2} \\ \hline
	: & : \\ \cline{2-2} 
	: & {\tiny 4n-q} \\ \cline{2-2}
	: \\ \cline{1-1}
	{\tiny 2n} \\ \cline{1-1} \end{tabular} \quad 
$\tensor$
\quad \begin{tabular}{|c|c|c|c|c|c|} \hline
	{\tiny 1} & {\tiny 2} & \multicolumn{3}{l}{$\dotsi$} &\vline\; {\tiny 2n}\\ \hline
	{\tiny 2n+1}& {\tiny 2n+2}& .. & {\tiny 4n-q} \\ \cline{1-4} \end{tabular}
\quad
$\too$
\quad \begin{tabular}{|c|c|} \hline
	{\tiny 1}  &{\tiny 2n+1}  \\ \hline
         {\tiny 2} & {\tiny 2n+2} \\ \hline
        : & : \\ \cline{2-2}
        : & {\tiny 4n-q} \\ \cline{2-2}
	: & 	{\tiny 1}  \\ \cline{2-2}
	: & : \\ \hline
	{\tiny 2n} & {\tiny q} \\ \hline \end{tabular}
\quad $\tensor$ \quad
	\begin{tabular}{|c|c|c|c|c|c|c|} \hline
	{\tiny 1} & {\tiny 2} & \multicolumn{4}{c}{$\dotsi$}\vline &{\tiny 2n}  \\ \hline
	{\tiny 2n+1} &{\tiny 2n+2}  & .. & {\tiny 4n-q}& {\tiny 1}& .. & {\tiny q} \\ \hline
	\end{tabular}

\vskip .2in

Note that $\gamma\gamma=1$.

Let $V$ be an orthogonal representation
of the fundamental group $\pi_1(M)$. Then we can consider differential forms 
on $M$ with values in $V$.
Then we can define the quaternionic analytic torsion using this complex
and the Laplacian $\Delta_q^D=D\gamma d\gamma + \gamma d\gamma D$. 
Denoting $\gamma d\gamma =\vardelta$, we can write
$\Delta_q^D = D\vardelta + \vardelta D$.

\begin{Def}
If $M$ is a quaternionic manifold and $V$ is a flat vector bundle,
then the quaternionic analytic torsion is defined by
$$\tau_\bH(M,V)=\prod_{q\ge 0}({{\text {det}}'\Delta_q^D})^{(-1)^{q+1}q/2},$$
where $\det'\Delta_q^D$ is the regularized determinant of $\Delta_q^D$.
\end{Def}

\subsection{Invariance of the Quaternionic Analytic Torsion}\label{varq}

Let $M^{4n}$ be a $4n$-dimensional quaternionic manifold and $g_u$ be a 
one parameter family of metrics on $M$ compatible to the almost quaternionic 
structure parametrized by $u$. We denote $\dot \gamma  =\frac {\partial}{\partial u}
\gamma $ and $\alpha=\gamma ^{-1}\dot\gamma =\gamma \dot\gamma =-\dot\gamma 
\gamma $. For simplicity we denote $\Delta_q^D$ as $\Delta_q$.
Then 
\begin{align*}
\dot\Laplace&=D\dot\gamma d\gamma +D\gamma d\dot\gamma +\dot\gamma d\gamma D+\gamma d\dot\gamma D\\
            &=D\dot\gamma \gamma \gamma d\gamma +D\gamma d\gamma \gamma \dot\gamma +\dot\gamma \gamma \gamma d\gamma D+\gamma d\gamma \gamma \dot\gamma D\\
            &=-D\alpha\vardelta+D\vardelta\alpha-\alpha\vardelta D
+\vardelta\alpha D
\end{align*}

With repeated application of $\tr(AB)=\tr(BA)$ for bounded operators 
as in the Ray and Singer's work,
we have
$$\tr(e^{t\Delta_q}\alpha\delta D)=\tr(e^{\frac12t\Delta_q}\alpha\delta D
e^{\frac12t\Delta_q})=\tr(\delta D d^{t\Delta_q \alpha}),$$
$$\tr(e^{t\Delta_q}\delta\alpha D)=\tr(D e^{t\Delta_q} \delta\alpha),\quad$$
$$ \tr(e^{t\Delta_q}D\alpha\delta)=\tr(\delta e^{t\Delta_q}D \alpha).$$

Hence we have
\begin{align*}
\frac{d}{du}\tr(e^{t\Laplace_u})&=t \tr(\dot\Laplace e^
{t\Laplace_u})\\
&=t \tr(-\alpha\vardelta D e^{t\Laplace_q}+
\alpha D\vardelta e^{t\Laplace_{q+1}}-\alpha\vardelta D e^{t\Laplace_{q-1}}
+\alpha D \vardelta e^{t\Laplace_q}).
\end{align*}

Hence,
\begin{align*}
\sum_{q=0}^{2n}(-1)^q q \tr(\dot\Laplace e^{t\Laplace_q})&=
\sum^{2n}_{q=0}(-1)^q\tr(\alpha \Laplace_q e^{t\Laplace_q})\\
&=\sum^{2n}_{q=0}(-1)^q\frac{d}{dt}\tr(\alpha e^{t\Laplace_q})
\end{align*}

Let $\zeta_1$ and $\zeta_2$ denote two zeta functions corresponding to 
two representations $V_1$ and $V_2$, and $\Delta^{(1)}$ and $\Delta^{(2)}$ 
denote corresponding Laplacians. Then the variation of the ratio between 
two quaternionic analytic torsions is given by

\begin{align*}
\frac{d}{du}\log \tau_\bH(M,V_1)/\tau_\bH(M,V_2)
&=\frac12\sum^{2n}_{q=0}\frac {d}{du}(-1)^qq
\left(\zeta_{1q}'(0)-\zeta_{2q}'(0)\right)\\
&=\frac{d}{ds}|_{s=0}\frac12\frac{1}{\Gamma(s)}\frac{d}{du}\sum_{q=0}^{2n}
(-1)^qq\int_0^{\infty}t^s\tr(e^{t\Delta_u})dt\\
&=\frac{d}{ds}|_{s=0}\frac12\frac{1}{\Gamma(s)}\sum^{2n}_{q=0}(-1)^q
\int^\infty_0 \frac{d}{dt}\tr(\alpha e^{t\Laplace_{(1)}}-\alpha
e^{t\Laplace_{(2)}})dt\\
&=\tr(\alpha(P^{(1)}-P^{(2)}))
\end{align*}
where $P$ denotes the orthogonal projection onto the space of harmonic 
forms corresponding to the Laplacian of proper degree.

Hence if the space of the harmonic forms is trivial, i.e., if the cohomology
$H^q_A(M,V_1)$ and $H^q_A(M,V_2)$ of the quaternionic complexes are
all trivial, then the variation is 0.

Here, we have used the fact that
$$\tr(e^{t\Laplace_{(1)}}-\alpha
e^{t\Laplace_{(2)}})=O(e^{-c/t}) \quad \text {as} \quad t \rightarrow 0$$
and that tr$(\alpha e^{t\Laplace_{(1)}}-\alpha
e^{t\Laplace_{(2)}})$ decreases exponentially since $\Laplace$ is strictly 
negative, hence the integration in the second line in the above equation
can be differentiated inside the integral sign.

Hence we have proven
\begin{Thm}
For a quaternionic manifold $M$ 
the ratio 
$$\tau_\bH(M,V_1)/\tau_\bH(M,V_2)$$
of two quaternionic analytic torsions corresponding to unitary representations
$V_1$ and $V_2$ is independent of the choice of the compatible metric if  the cohomology group $H^*_A$
is trivial.
\end{Thm}

\begin{Rem}
Recall that the Riemannian structure on $M^4$ determines a quaternionic
structure and the quaternionic complex for $M$ is the same as the
self-dual complex.
It is not difficult to see from the inclusion
\begin{equation*}
\begin{tabular}{ccc}
$SO(4)$ & $\subset$ & $GL(1,\bH)Sp(1)$\\
$\parallel$ & & $\parallel$\\
$Sp(1)Sp(1)$ & $\subset$ & $\bH^\times\cdot Sp(1)$
\end{tabular}
\end{equation*}
that two Riemannian metrics $g$ and $g^{\prime }$
determines the same quaternionic structure on $M$ provided that they are
conformal to each other, that means $g^{\prime }=e^{2v}g$ for some smooth
function $v$ on $M.$ 

As a result, the regularized
determinant of the self-dual complex of  $M$ depends only on the conformal
class of the metric. In the next section, we shall see that the self-dual
analytic torsion is always a conformal invariant quantity.
\end{Rem}

Recall that the real analytic torsion vanishes on complex manifolds and the
complex analytic torsion vanishes on HyperK\"{a}hler manifolds. It is
natural to ask if the quaternionic analytic torsion vanishes for certain special
class of HyperK\"{a}hler manifolds. In fact to look for the corresponding
vanishing result, we should consider analytic torsions defined for manifolds
admitting {\it{commuting}} complex structures instead of anti-commuting
complex structures as in HyperK\"ahler manifolds.  Manifolds admitting commuting complex structures are studied in recent
years by some string theorists (For example, see Rocek \cite{R}) as
candidates for mirrors of rigid Calabi-Yau manifolds of dimension three.
Examples of such manifolds includes the Hopf surface or products of the Hopf
surface with any complex manifold. (Recall that the Hopf surface is $
\bC^{2*}/{\left( z_1,z_2\right) \thicksim
\left( 2z_1,2z_2\right) }$ with complex structure descended from $\bC^2.$%
)

\section{Self Dual Analytic Torsions}

As we discussed before, the quaternionic analytic torsion for 4-dimensional
quaternionic manifold is the same as the determinant of the self-dual 
complex of the manifold.

In this section, we shall define the self-dual analytic torsion for 4n-dimensional manifolds and show that the ratio between two analytic torsions
$\tau _{\text{SD}}\left( M,V_1,g\right) /$ 
$\tau _{\text{SD}}\left( M,V_2,g\right) $ corresponding to two flat bundles
$V_1$ and $V_2$ is a conformally invariant quantity.
Then we study the K\"ahler surface in more details and compare various torsions
on them.

\subsection{Definition and Invariance of Self-dual Analytic Torsion}
Let $M$ be a $4n$-dimensional Riemannian manifold with metric $g$ and $V$ be 
the flat vector bundle
associated with an orthogonal representation of $\pi_1(M)$.
Then we have
\begin{Def}
The self-dual analytic torsion is defined by
$$\log\tau_{\text{SD}}(M,V,g)=\sum\limits_{q=0}^{2n-1}\left(
-1\right) ^{q+1}q\log \det \left(\bigtriangleup _q\right) -n\log
\det \left(\bigtriangleup _{2n}\right).$$ 
where $\Delta_q$ is the Laplacian on $q$-forms on $M$ with values in $V$.
\end{Def}
\begin{Prop}
The self-dual analytic torsion is the regularized determinant of the following
self-dual complex
\[
0\rightarrow \Omega ^0(M,V) \stackrel{d}{\rightarrow }\dots\stackrel{d}{\rightarrow 
}\Omega ^{2n-1}(M,V) \stackrel{\sqrt{2}P_{+}d}{\rightarrow }%
\Omega _{+}^{2n}(M,V) \stackrel{}{\rightarrow }0
\]
where $P_+=(*+id)/2$ is the projection to the space of self-dual forms.
\end{Prop}

This is an elliptic complex which can be regarded as half of the deRham
complex. The other half is the anti-self-dual complex,
\[
0\rightarrow \Omega _{-}^{2n}(M,V) \stackrel{\ }{\stackrel{d}
{\rightarrow}\dots}\stackrel{d}{\rightarrow }\Omega ^{4n-1}(
M,V) \stackrel{d}{\rightarrow }\Omega ^{4n}(M,V)
\rightarrow 0.
\]
The indexes of these two complexes are respectively
$\frac12(\chi(M,V)+sign(M,V)) $ 
and $\frac12( \chi(M,V) -sign(M,V))$. Their sum is just
the Euler characteristic of $M$ with twisted coefficient $V$.
In a similiar vein, we can define the {\it{anti-self-dual analytic torsion}}
 $\tau_{\text{ASD}}(M,V)$. By the vanishing result of Ray and Singer, we get 
$$ \tau_{\text{SD}}(M,V)\cdot\tau_{\text{ASD}}(M,V)=1.$$
Hence the anti-self-dual analytic torsion doesn't give any new information.

\begin{proof}(of proposition)
If we denote $\left( \sqrt{2}P_{+}\circ d\right) \circ \left( \sqrt{2}%
P_{+}\circ d\right) ^{*}$ by $\bigtriangleup _{+},$ then it coincides with
the ordinary Laplacian $\bigtriangleup _{2n}$ on self-dual forms. It is 
because
\begin{align*}
\bigtriangleup _{2n}\phi & =  -d*d*\phi -*d*d\phi \\ 
& =  -\left( *+1\right) d*d\phi \\ 
& =  -2P_{+}d*d\phi
\end{align*}
where $\phi \in \Omega _{+}^{2n}\left( M\right)$ and
on the other hand, we have $(P_{+}\circ d)^{*}\phi =\delta \phi
=-*d\phi $ because 
\begin{align*}
\int \left\langle \psi ,\left( P_{+}\circ d\right) ^{*}\phi \right\rangle & =
 \int \left\langle \left( P_{+}\circ d\right) \psi ,\phi \right\rangle  
 =  \int \left\langle d\psi ,\phi \right\rangle  
 =  \int \left\langle \psi ,\delta \phi \right\rangle.
\end{align*}
Hence
\begin{align*}
\left( P_{+}\circ d\right) \circ \left( P_{+}\circ d\right) ^{*}\phi  =  
-P_{+}\circ d*d\phi =  \frac 12\bigtriangleup _{2n}\phi
\end{align*}
and we have 
\begin{align*}
\tr\,e^{t\bigtriangleup _{+}}|_{\Omega _{+}^{2n}(M) } & =  
\tr\,e^{t\bigtriangleup _{2n}}|_{\Omega _{+}^{2n}(M) } \\ 
\  & =  \tr\,e^{t\bigtriangleup _{2n}}\circ P_{+}|_{\Omega ^{2n}(
M) } \\ 
& =  \tr\,e^{t\bigtriangleup _{2n}}\circ (*+1)/2 |_{\Omega
^{2n}( M) } \\ 
& =  \frac 12\tr\,e^{t\bigtriangleup_{2n}}|_{\Omega ^{2n}(M)
}+\frac 12\tr\,e^{t\bigtriangleup _{2n}}*|_{\Omega ^{2n}(M) }\ 
\end{align*}
One can show that the
variation of $\tr\,e^{t\bigtriangleup _{2n}}*|_{\Omega
^{2n}(M,V) }$ under a conformal change of metrics is zero and hence 
a conformally invariant. In fact, it is a homotopy invariant, namely the 
signature of $M$ with coefficients in $V$: 
>From McKean Singer's formula, we have $sign(M,V)
=\sum\limits_{q=0}^{4n}\tr e^{t\bigtriangleup _q}*|_{\Omega ^q\left(
M,V\right) }.$ However, all non-middle dimensional contributions to the
summation cancel each other and gives $sign\left( M,V\right) =$ $%
\tr\,e^{t\bigtriangleup _{2n}}*|_{\Omega ^{2n}\left( M,V\right) }.$
Alternatively, 
\[
\left[ d,\delta \right] :\Omega _{+}^{2n}\left( M,V\right) \rightarrow
\Omega _{-}^{2n}\left( M,V\right) 
\]
is an isometry which preserves eigenspaces of Laplacian with non-zero
eigenvalue. Therefore, $\tr\,e^{t\bigtriangleup _{2n}}*|_{\Omega ^{2n}\left(
M,V\right) }=\tr*|_{H^{2n}\left( M,V\right) }$ which is clearly the twisted
signature of $M$ with coefficients in $V$.

However, for the regularized determinant of an operator, we discard all the
zero eigenvalues. So $\tr\,e^{t\bigtriangleup _{2n}}*|_{\Omega ^{2n}\left(
M,V\right) }$ does not contribute to the torsion. Hence we 
can see that the analytic torsion defined by the elliptic complex 
\begin{align*}
0\rightarrow \Omega^0( M,V) \stackrel{d}{\longrightarrow}\dots
\stackrel{d}{\longrightarrow}\Omega^{2n-1}(M,V) 
\stackrel{\sqrt{2}P_{+}d}{\longrightarrow}\Omega_{+}^{2n}
(M,V) \rightarrow 0
\end{align*}
is given by
\[
\log\tau_{\text{SD}}(M,V)=\sum\limits_{q=0}^{2n-1}\left(
-1\right) ^{q+1}q\log \det \left(\bigtriangleup _q\right) -n\log
\det \left(\bigtriangleup _{2n}\right) 
\]
\end{proof}
Next we show that the self-dual analytic torsion is a conformal invariant quantity.
\begin{Thm}
For any two orthogonal flat bundles $V_1$ and $V\ _2$ over $M$ with trivial
cohomology group. The ratio of the two self-dual analytic torsion $\tau _{\text{SD}}\left( M,V_1,g\right) /$ $\tau _{\text{SD}}\left(
M,V_2,g\right) $ depends only on the conformal class of the metric $g$.
\end{Thm}

\begin{proof}
To compute the variation of the self-dual analytic torsion under conformal change
of metrics $g\left( u\right) $, we need to know $\frac \partial {\partial
u}\tr\,e^{t\bigtriangleup _q\left( u\right) }.$ The following formulae (from
Ray-Singer's paper) are valid for any change of Riemannian metric $\frac{%
\partial g}{\partial u}$%
\[
\frac \partial {\partial u}\tr\,e^{t\bigtriangleup _q\left( u\right)
}=t\tr\left( \,e^{t\bigtriangleup _q\left( u\right) }\dot{\bigtriangleup}%
_q\right) 
\]
where $\dot{\bigtriangleup}_q=\frac{d\bigtriangleup _q\left( u\right) }{du}%
=\alpha \delta d-\delta \alpha d+d\alpha \delta -d\delta \alpha $ and $%
\alpha =*^{-1}\dot{*}$ $=*^{-1}\frac{\partial *}{\partial u}.$ Therefore, 
\[
\tr\,\,e^{t\bigtriangleup _q}\dot{\bigtriangleup}_q=\tr\,\,e^{t\bigtriangleup
_q}\delta d\alpha +\tr\,\,e^{t\bigtriangleup _{q-1}}\delta d\alpha
-\tr\,\,e^{t\bigtriangleup _q}d\delta \alpha -\tr\,\,e^{t\bigtriangleup
_{q+1}}d\delta \alpha 
\]
because $d\bigtriangleup _q=\bigtriangleup _{q+1}d$ and $\delta
\bigtriangleup _q=\bigtriangleup _{q-1}\delta .$

>From Ray-Singer's paper, we have

\begin{align*}\label{lll}
&\sum\limits_{q=0}^k\left( -1\right) ^qq\tr e^{t\bigtriangleup _q}\dot{ \bigtriangleup}_q   \\ 
&=\sum\limits_{q=0}^k(-1) ^{q+1}\tr\,e^{t\bigtriangleup_q}\bigtriangleup _q\alpha+\left( -1\right) ^k\left( k+1\right) \tr\,e^{t\bigtriangleup _k}\delta
d\alpha +\left( -1\right) ^{k+1}k\tr\,e^{t\bigtriangleup _{k+1}}d\delta
\alpha  
\end{align*}
for any integer $k.$ In their paper, they study the case when $k$ is the
dimension of $M$ and define a differentiable invariant. For our purposes, we
shall take $k=2n-1$ when $\dim M=4n.$ To counter the error term 
$-2n\tr\,e^{t\bigtriangleup _{2n-1}}\delta d\alpha +\left( 2n-1\right)
\tr\,e^{t\bigtriangleup_{2n}}d\delta \alpha $ we need to consider 
$$n\cdot \tr\,\,e^{t\bigtriangleup_2n}\dot{\bigtriangleup}_2n+\sum
\limits_{q=0}^{2n-1}\left( -1\right) ^qq\tr\,\,e^{t\bigtriangleup _q}
\dot{\bigtriangleup}_q.$$ 
Now 
\[
\tr\,e^{t\bigtriangleup_{2n}}\dot{\bigtriangleup}_2n=\tr\,\,e^{t\bigtriangleup
_{2n-1}}\delta d\alpha -\tr\,\,e^{t\bigtriangleup _{2n+1}}d\delta \alpha 
\]
because $\dot{\alpha}=0$ on $\Omega ^{2n}\left( M\right) ,$ under conformal
change of metrics.

By applying above formulae, we get 
\begin{align*}
&n\tr\,e^{t\bigtriangleup_{2n}}\dot{\bigtriangleup}_{2n}
+\sum\limits_{q=0}^{2n-1}( -1) ^qq\tr\,e^{t\bigtriangleup_q}
\dot{\bigtriangleup}_q \\ 
&=\sum\limits_{q=0}^{2n-1}(-1)^{q+1}\tr\,e^{t\bigtriangleup
_q}\bigtriangleup_q\alpha -2n\tr\,e^{t\bigtriangleup_{2n-1}}\delta
d\alpha +n( \tr\,\,e^{t\bigtriangleup _{2n-1}}\delta d\alpha
-\tr\,\,e^{t\bigtriangleup _{2n+1}}d\delta \alpha )  \\ 
&=  \sum\limits_{q=0}^{2n-1}( -1) ^{q+1}\tr\,e^{t\bigtriangleup
_q}\bigtriangleup _q\alpha - n( \tr\,\,e^{t\bigtriangleup
_{2n-1}}\delta d\alpha +\tr\,\,e^{t\bigtriangleup _{2n+1}}d\delta \alpha) 
\end{align*}

Since $\tr\,\,e^{t\bigtriangleup_{2n-1}}\delta d\alpha+\tr\,\,e^{t\bigtriangleup
_{2n+1}}d\delta \alpha =0$, we have
$$n \tr\,\,e^{t\bigtriangleup_{2n}}\dot{\bigtriangleup}
_{2n}+\sum\limits_{q=0}^{2n-1}\left( -1\right) ^qq\tr\,\,e^{t\bigtriangleup _q}
\dot{\bigtriangleup}_q=\sum\limits_{q=0}^{2n-1}\left( -1\right)
^{q+1}\tr\,e^{t\bigtriangleup _q}\bigtriangleup _q\alpha.$$

Therefore, for $s$ with sufficiently large $\func{Re}s$, we have
\begin{align*}
&\frac d{du}\left[ \sum\limits_{q=0}^{2n-1}(-1)^qq\zeta
_q(s)+n\zeta_{2n}(s)\right]  \\
& =  \frac 1{\Gamma \left( s\right) }\int\limits_0^\infty t^s
\left( n\tr\,\,e^{t\bigtriangleup _2n}
\dot{ \bigtriangleup}_{2n}+\sum\limits_{q=0}^{2n-1}\left( -1\right)
^qq\tr\,\,e^{t\bigtriangleup _q}\dot{\bigtriangleup}_q\right) \ \,dt \\ 
& =  \sum\limits_{q=0}^{2n-1}(-1) ^{q+1}\frac 1{\Gamma(s)}\int\limits_0^\infty t^s\tr\,e^{t\bigtriangleup_q}
\bigtriangleup _q\alpha dt \\ 
& =  \sum\limits_{q=0}^{2n-1}(-1)^{q+1}\frac 1{\Gamma 
(s)}\int\limits_0^\infty t^s
\frac d{dt}\tr(\,e^{t\bigtriangleup _q}\alpha) dt
\end{align*}

>From \cite{RS2} it follows that 
\[
\frac 1{\Gamma(s) }\int\limits_0^\infty t^s\frac
d{dt}\tr\,\left( e^{t\bigtriangleup _q^{(1) }}-e^{t\bigtriangleup
_q^{(2) }}\right) \alpha dt
\]
defines an analytic function of $s$ which vanishes at $s=0$ where $%
\bigtriangleup _q^{(i)}$ denotes the Laplacian on $M$ with
coefficient in $V_i.$ We have 
\begin{align*}
\frac d{ds}|_{s=0}\left[ \frac 1{\Gamma(s)}\int\limits_0^\infty 
t^s\frac d{dt}\tr\,\left( e^{t\bigtriangleup _q^{(1)}}
-e^{t\bigtriangleup _q^{(2) }}\right) \alpha dt\right]   
=  \tr(P^{(1) }-P^{(2)}) \alpha 
\end{align*}
where $P^{(i)}$ is the harmonic projection operator for $i=0,1.$
>From our assumption, $P^{(i) }=0$ and hence we have shown that 
\[
\frac d{du}\sum\limits_{i=1}^2(-1) ^i\left(
\sum\limits_{q=0}^{2n-1}(-1) ^qq\zeta_q^{^{\prime
}}(0,V_i)+n\zeta_{2n}^{^{\prime }}(0,V_i)\right) =0
\]
That is $\tau _{\text{SD}}\left( M,V_1,g\right) /$ $\tau _{\text{SD}}\left(
M,V_2,g\right) $ is independent of $u.$ Hence, we have our theorem.
\end{proof}

\begin{Rem}
Since $\alpha =0$ on middle dimensional forms, we get
$$\tr|_{\Omega ^{2n}\left(
M,V\right) }\left( P^{(1) }-P^{(2) }\right) \alpha =0$$ 
Therefore, it is enough for us to assume that $H^k\left( M,V_i\right) $
are zero for $k<2n$ for the result to hold.
\end{Rem}

\begin{Rem}
Notice that if a closed $4n$-dimensional manifold $M$ is a complex
manifold, then the complexification of the self-dual complex has a
subcomplex, namely the Dolbeault complex:
\[
\begin{array}{cccccc}
0\rightarrow&\negthickspace\Omega ^0(M)_{\bC}\stackrel{d}{\rightarrow }
& 
\negthickspace\dots\rightarrow  & \negthickspace\Omega ^{2n-1}(M) _{\Bbb{C}}\stackrel{\sqrt{%
2}P_{+}d}{\rightarrow } & \negthickspace\Omega _{+}^{2n}(M)_{\Bbb{C}%
}\rightarrow  & 0 \\
& \parallel   &  & \cup  & \cup  &  \\
0\rightarrow&\negthickspace\Omega ^{0,0}(M) \stackrel{\bar{\partial}}{%
\rightarrow } &  \negthickspace\dots\rightarrow  & \negthickspace\Omega ^{0,2n-1}(M)
\stackrel{\bar{\partial}}{\rightarrow } & \negthickspace\Omega ^{0,2n}(M)
\rightarrow  & 0.
\end{array}
\]

All these inclusions are trivial except possibly the last one which again
can be checked directly (without assuming K\"{a}hlerian property of $M$).
\end{Rem}

\begin{Rem}
Similiarly if a closed $4n$-dimensional manifold $M$ is a
quaternionic manifold, then the complexification of the self-dual complex
again has a subcomplex, namely the quaternionic complex:
\end{Rem}
$$\small{
\begin{array}{cccccc}
0\rightarrow&\negthickspace\Omega ^0(M) _{\Bbb{C}}\stackrel{d}{\rightarrow }
&
\negthickspace\dots\rightarrow  & \negthickspace\Omega ^{2n-1}(M) _{\Bbb{C}}\stackrel{\sqrt{%
2}P_{+}d}{\rightarrow } &\negthickspace \Omega _{+}^{2n}(M) _{\Bbb{C}%
}\rightarrow  & 0 \\
&\parallel  &  & \cup  & \cup  &  \\
0\rightarrow&A^0(M) \stackrel{D}{\rightarrow } & \dots\rightarrow  & A^{2n-1}(M) \stackrel{D}{\rightarrow } & A^{2n}(M) \stackrel{D}{%
\rightarrow } & 0.
\end{array}
}
$$
\section{Acknowledgement}

{ The first author would like to thank I. M. Singer.}


\end{document}